\documentclass[10pt]{article}
\usepackage{latexsym}
\usepackage{amsfonts}
\usepackage{graphicx}
\title{Whitehead method and Genetic Algorithms}
\author{Alexei D. Miasnikov, Alexei G. Myasnikov}
\date{}
\newtheorem{theorem}{Theorem}

\newtheorem{conjecture}{Conjecture}
\newtheorem{problem}{Problem}

\newtheorem{conclusion}{Conclusion}
\newtheorem{remark}{Remark}

\begin{document}
\newcommand{\PA}{A }
\newcommand{\PB}{B }

 \pagestyle{myheadings} \markright{{\scriptsize  \textsf{
   Whitehead Method and Genetic Algorithms,
    A. D. Miasnikov,    A. G. Myasnikov}  $\bullet$ 04/17/2003}}

\maketitle

\begin{abstract}
In this paper we discuss a genetic version (GWA) of the
Whitehead's algorithm, which is one of the basic algorithms in
combinatorial group theory.  It turns out that GWA is surprisingly
fast and outperforms  the standard Whitehead's algorithm  in free
groups of rank $ \geq 5$. Experimenting with GWA we  collected an
interesting numerical data that clarifies the time-complexity of
the Whitehead's Problem in general. These experiments led us  to
several mathematical conjectures. If confirmed they will shed
light on hidden mechanisms of Whitehead Method and geometry of
automorphic orbits in free groups.
\end{abstract}

\section{Introduction}

Genetic Algorithms have been introduced by J.H.Holland in 1975
\cite{Hol}. Since then they have been successfully applied in
solving  a number of numerical and combinatorial problems. In most
cases genetic algorithms are used in optimization problems
 when  searching for an optimal solution or its approximation
  (see, for example, survey \cite{Reeves}).

The first applications of genetic algorithms to abstract algebra
appeared in \cite{Me} and  \cite{MeAlex}, where we made some
initial attempts to  study Andrews-Curtis conjecture from
computational view-point. In the present paper we discuss a
genetic version of Whitehead algorithm, which is one of the basic
algorithms in combinatorial group theory.  It turns out that this
Genetic Whitehead Algorithm (GWA) is surprisingly fast and
outperforms  the standard Whitehead algorithm  in free groups of
rank $ \geq 5$. Experimenting with GWA we were able to collect an
interesting numerical data which clarifies the time-complexity of
Whitehead Problem in general. These experiments led us  to several
mathematical conjectures which we stated at the end of the paper.
If confirmed they will shed  light on hidden mechanisms of
Whitehead Method and geometry of automorphic orbits in free
groups. Actually, the remarkable performance of GWA has initiated
already investigation of automorphic orbits in free groups of rank
2 \cite{MS}, \cite{Khan}. Some of the conclusions that one can
draw from our experiments are worth to be mentioned here.

One unexpected outcome of our experiments is that the time
complexity functions of Whitehead's algorithms in all their
variations does not depend "essentially" on the length of the
input words. We introduce a new type of size function (Whitehead's
Complexity function) on input words which allows one to measure
adequately the time complexity of Whiteheads algorithms. This type
of size functions is interesting in its own right, it makes
possible to compare a given algorithm from a class of algorithms
${\cal K}$ with the best possible non-deterministic algorithm in
${\cal K}$.

This Whitehead's complexity function takes care of the observed
phenomena that most of the words in a given free group are already
Whitehead's minimal (have minimal length in their automorphic
orbit). Such words have Whitehead's complexity 0 and the
Whitehead's descent algorithm is meaningless for such words.

Another conclusion we made is that the actual generic (or average)
time complexity of the Whitehead's descent algorithm (on
non-minimal inputs, of course) is much less than of the standard
Whitehead's algorithm. Moreover, it does not depend on the rank
$r$ of the ambient free group $F_r$ exponentially, though the
standard one does. We believe that there exists a finite subset
$T_r$ (of polynomial size in $r$) of elementary Whitehead's
automorphisms in $F_r$ for which the classical Whitehead's descent
method does nor encounter any "picks" on the most inputs.

 Genetic Whitehead Algorithms GWA was designed and implemented
in 1999 and soon after some interesting facts transpired from
experiments. But only recently an adequate group-theoretic
language (average case complexity, generic elements, asymptotic
probabilities on infinite groups) was developed which would allow
one to describe the group-theoretic part of the observed
phenomena. We refer to papers \cite{BMS}, \cite{BMR},
\cite{KMSS1}, \cite{KMSS2} for details. On the other hand,  a
rigorous theory of genetic algorithms is not developed yet up to
the level which would explain fast performance of such heuristic
algorithms as GWA. In fact, we believe that thorough investigation
of particular genetic algorithms in abstract algebra might provide
insight to a general theory of genetic algorithms.

\section{Whitehead's method}

\subsection{Whitehead Theorem} \label{sub-sec:WTh}

Let $X = \{x_1, \ldots, x_n\}$ be a finite  set and $F = F_n(X)$ be a
free group with a basis $X$.  Put  $X^{ \pm 1} = \{ x^{\pm 1} | x \in
X\}$. We will represent elements of $F$ by reduced words in the
alphabet  $X^{\pm 1}$ (i.e., words without subwords $x x^{-1}, x^{-1}
x$ for any $x \in X$). For a word $u$ by $|u|$ we denote the length of
$u$, similarly, for a tuple $U = (u_1, \ldots, u_k) \in F^k$ we denote
by $|U|$ the {\it total length} $|U| = |u_1| + \ldots + |u_k|$.

 For an automorphism  $\varphi$ of $F$, and    $k$-tuples
 $U=(u_1,...,u_k), V=(v_1,...,v_k) \in F^k$ we  write
$U\varphi=V$ if $u_i \varphi = v_i$, $i=1,...,k$.

In 1936  J.H.C. Whitehead  introduced the following algorithmic
problem, which became a central problem of the theory of
automorphisms of free groups  \cite{W}.

\medskip
{\bf Problem  W}  {\it Given two tuples $U, V \in F^k$  find out if
there is an automorphism $\varphi \in Aut(F)$ such that $U\varphi =
V$.}

\medskip
In the same paper he showed (using a topological argument) that
this problem  can be solved algorithmically  and suggested an
algorithm to find such an automorphism $\varphi$ (if it exists).
To explain this method we need the following definition.
 An automorphism $t \in Aut(F)$  is called a  \emph{Whitehead automorphism}
 if it has  one of the following types:
\begin{enumerate}
\item [1)] $t$ permutes elements in $X^{\pm 1}$;
\item [2)] $t$ takes each element $x \in X^{\pm 1}$ to one of the elements
$x$, $xa$, $a^{-1} x$, or $a^{-1} x a$, where $x \neq a^{\pm 1}$ and $a \in
X^{\pm 1}$ is a fixed element.
\end{enumerate}

Denote by $\Omega_n = \Omega(F)$ the set of all Whitehead automorphisms
of a given free group $F = F_n(X)$. It follows from a result of Nielsen
that $\Omega_n$ generates $Aut(F_n(X))$ \cite{Nielsen}.

Let $T$  be a subset of $Aut(F)$. We say that tuples $U, V \in F^k$ are
$T$-equivalent ($U \sim_T V$) if there exists a finite sequence $t_1,
\ldots, t_m$ (where $t_i \in T^{\pm 1}$) such that $Ut_1 \ldots t_m =
V$. The $T$-equivalence class of a tuple $U$ is called the $T$-orbit
$Orb_T(U)$ of $U$.  If $T$ generates $Aut(F_n)$  then the equivalence
class of a tuple $U$ is called the {\it orbit} $Orb(U)$ of $U$. Now
Problem W can be stated as a membership problem for a given orbit
$Orb(U)$. By $U_{min}$ we denote any tuple of minimal total length in
the orbit $Orb(U)$, and by $Orb_{min}(U)$ - the set of all minimal
tuples $U_{min}$.

It is  convenient  sometimes to look at Whitehead problem from
graph-theoretic view-point. Denote by $\Gamma(F,k,T)$  the following
directed labelled graph: $F^k$ is the vertex set of $\Gamma$; two
vertices $U, V \in F^k$ are connected by a directed  edge from $U$ to
$V$ with label $t \in T$ if and only if $Ut = V$.   We refer to
$\Gamma_k(F) = \Gamma(F,k,\Omega)$ as to {\it Whitehead graph} of $F$.
In the case when $k = 1$ we write $\Gamma(F)$ instead of $\Gamma_1(F)$.
 Obviously, $V \in Orb(U)$ if and only if $U$
and $V$ are in the same connected component of $\Gamma_k(F)$.


The following theorem is one of the fundamental results in
combinatorial group theory.

\begin{theorem}
\label{T1}
(Whitehead \cite{W}). Let $U,V \in F_n(X)^k$ and $V \in Orb(U)$. Then: \\
1) if $|U|>|V|$, then there exists $t \in \Omega_n$ such that
\[|U|>|U
t|;\] 2) if $|U| = |V|$, then there exist $t_1, \ldots, t_m \in
\Omega_n$ such that
\[U t_1 ... t_m = V\]
and $|U|=|U t_1|=|U t_1 t_2| = ... = |U t_1 t_2 ... t_m| = |V|.$
\end{theorem}

 In view of  Theorem \ref{T1}  Problem W  can be
divided into two subproblems:

\medskip {\bf Problem A} {\it  For a tuple $U \in F^k$ find a sequence
$t_1, \ldots, t_m \in \Omega_n$ such that $Ut_1 \ldots t_m = U_{min}$.
}

\medskip
{\bf Problem B} {\it  For tuples  $U, V \in F^k$ with
$$|U| = |U_{min}| = |V_{min}| = |V|$$
find a sequence $t_1, \ldots, t_m \in \Omega_n$ such that $Ut_1 \ldots
t_m = V$.}

\medskip
Theorem \ref{T1} gives a solution to the both problems above, and
hence to  Problem W.

\subsection{Whitehead Algorithm}
\label{sub-sec:WA}

The procedures  described  below give algorithmic solutions to the
Problems A and B, together they are known as {\it Whitehead
Algorithm } or {\it Whitehead Method}.

\subsubsection{Decision algorithm for Problem  A}
\label{sub-sec:DWA}

Following Whitehead  we describe below a deterministic decision
algorithm for Problem A, we refer to this algorithm (and  to various
its modifications)  as to DWA. This algorithm executes consequently
 the following

 \medskip
 \noindent {\it "Elementary Length Reduction Routine"} (ELR):

\begin{quote}
 Let $U \in F^k$. ELR finds $t \in \Omega_n$ with $|Ut| < |U|$ (if it exists).
 Namely,  ELR performs the following search. For each $t
\in \Omega_n$ compute the length of the tuple $Ut$ until $|U| > |Ut|$,
then put $t_1 = t, U_1 = Ut_1$ and output $U_1$. Otherwise stop and
output  $U_{min} = U$.
\end{quote}

 DWA  performs  ELR on  $U$, then performs ELR
on $U_1$, and so on, until a minimal tuple $U_{min}$ is found. We
refer to   algorithms of this type as to  {\em  Whitehead's
descent method} with respect to the set $\Omega_n$.

Clearly, there could be at most $|U|$ repetitions of ELR:
\[|U| > |Ut_1|>...>|Ut_1 ... t_l| = U_{min}, \ \ \ l \leq |U|.\]
 The sequence $t_1, \ldots, t_l$ is a  solution to Problem A. Notice,
 that the iteration procedure above simulates the classical
 gradient descent method ($t_1$ is the best direction from $U$,
 $t_2$ is the best direction from $U_1$, and etc.).

\subsubsection{Decision algorithm for Problem  B.}
\label{sub-sec:DWB}

Here we describe a deterministic decision algorithm for Problem B,
which is also due to Whitehead. In the sequel we refer to this
algorithm (and its variations) as to DWB.

 Let $U, V \in F^k$. DWB  constructs  $Orb_{min}(U)$ (as well
as $Orb_{min}(V)$) by repeating consequently the following

\medskip\noindent
{\it "Local Search Routine"} (LS):
\begin{quote}
 Let $\Omega_n = \{t_1,
\ldots, t_m\}$ and $\Delta$ be a finite graph with vertices from
$F^k$. Given a vertex $W$ in $\Delta$ the local search at $W$
results in a graph $\Delta_W$ which contains $\Delta$. We define
$\Delta_W$ recursively. Put $\Gamma_{0} = \Delta$, and suppose
that  $\Gamma_{i}$ has been already constructed. If $|Ut_{i+1}| =
|U|$ and $Ut_{i+1}$ does not appear in $\Gamma_{i}$ then add
$Ut_{i+1}$ as a new vertex to $\Gamma_{i}$, also add a new edge
from $U$ to $Ut_{i+1}$ with label $t_{i+1}$, and denote the
resulting graph by $\Gamma_{i+1}$. Otherwise, put $\Gamma_{i+1} =
\Gamma_{i}$. The routine stops in $m$ steps and results in a graph
$ \Gamma_{m}$. Put $\Delta_W = \Gamma_{m}.$ \end{quote}

 The construction of $Orb_{min}(U)$  is a variation of the standard

 \medskip\noindent {\it  "Breadth-First Search Procedure"} (BFS):

\begin{quote}
 Start with a graph $\Delta_0$ consisting of a single vertex
$U$. Put $\Delta_1 = (\Delta_0)_W$ and "mark" the vertex $U$.  If a
graph $\Delta_i$ has been constructed, then take any unmarked vertex
$W$ in $\Delta_i$ within the {\it shortest} distance from $U$, put
$\Delta_{i+1}= (\Delta_i)_W$, and mark the vertex $W$.
\end{quote}

Since $Orb_{min}(U)$ is finite BFS terminates, say in $l$ steps, where
$$l \leq |Orb_{min}(U)||\Omega_n|$$
 It is easy to see that $\Delta_l$ is
a tree, containing all vertices from $Orb_{min}(U)$. This implies that
$V \in Orb_{min}(U)$ if and only if $V \in \Delta_l$. Moreover, the
unique path connecting $U$ and $V$ in $\Delta_l$ is a shortest path
between $U$ and $V$ in $Orb_{min}(U)$, and the sequence of labels along
this path is a  sequence of Whitehead automorphisms (required in
Problem B) that connects $U$ and $V$ inside $Orb_{min}(U)$.

From the computational view-point it is more efficient to start
building maximal trees in both graphs $Orb_{min}(U)$ and $Orb_{min}(V)$
simultaneously, until a common vertex occurs.

\begin{figure}
\centerline{\includegraphics{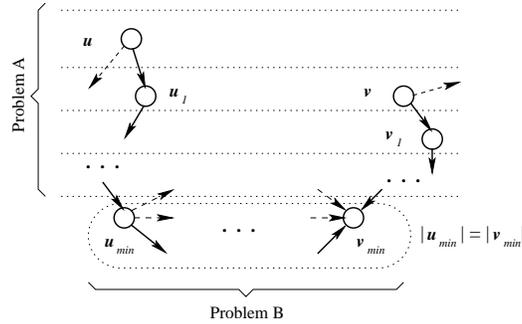}} \caption{Whitehead
Method.} \label{fig:whiteheadgraph}
\end{figure}

\subsection{Estimates for the time-complexity of the  Whitehead's
algorithms.} \label{sub-sec:est}

\subsubsection{Algorithm DWA.}
\label{sub-sec:comp-DWA}

It is easy to see that transformations of the type 1) cannot reduce the
total length of a tuple. Hence, to solve Problem A one needs only
Whitehead automorphisms of the type 2). It is not hard to show that
there are
 $$A_n = 2n4^{(n-1)} - 2n$$
  non-trivial Whitehead automorphisms of the type 2).


In the worst-case scenario to perform ELR it requires $A_n$ executions
of the following

\medskip\noindent {\it Substitution Routine} (SR):

\begin{quote}
For a given automorphism $t$ of the type 2) make a substitution $x
\rightarrow xt$ for each occurrence of each $x \in X^{\pm 1}$ in $U$,
and  then make all possible cancellations. \end{quote}

 Since the length of the word $xt$ is at most 3 the time needed to perform
 this routine is bounded from above by $c|U|$, where $c$ is a constant
 which does not depend on $|U|$ and the rank of $F$.
 Since DWA executes ELR at most $|U|$ times  the time-complexity function
  of DWA is bounded from above by
 $$cA_n|U|^2 = c(2n4^{n-1} - 2n)|U|^2,$$
 This bound depends exponentially on the rank
$n$ of the group $F = F_n(X)$.  For example, if $k = 1$,  $n = 10$, and
$|U| = 100$, the estimated number of steps for DWA  is bounded above by
 $$c(20 \cdot 4^9 - 20) 100^2 > c(5 \cdot 10^{10}).$$
 Whether this bound is tight in the
worst case is an open question.
In any event, computer experiments which we ran  on a dual Pentium
III, 700 Mhz processor computer with 1Gb memory show (see Table
\ref{tab:table1}) that the standard DWA  cannot find $U_{min}$ on
almost all inputs $U$ which are pseudo-randomly generated
primitive elements  of length more then 100 in the group $F_{10}$,
while working non-stop for more than an hour.

The accuracy of the bound depends on how many automorphisms from
$\Omega_n$ do reduce the length of a given input $U$. To this end,
put
$$LR(U) = \{t \in \Omega_n \mid |Ut| < |U|\}$$
Now, the number of steps that ELR performs on a {\em worst-case}
input $U$ is bounded from above by
$$\max \{A_n - |LR(U)|, 1\}$$
(if the ordering of $\Omega_n$ is such that all automorphisms from
$LR(U)$ are located at the end of the list $\Omega_n = \{t_1,
\ldots, t_m\}$).

If we assume that the automorphisms from $LR(U)$ are distributed
uniformly in the list $\Omega_n$ then DWA needs
$$ A_n^{\prime} = \frac{A_n}{|LR(U)|}$$
steps on average to find a length reducing automorphism for $U$.

\begin{figure}[ht]
\centerline{  \includegraphics{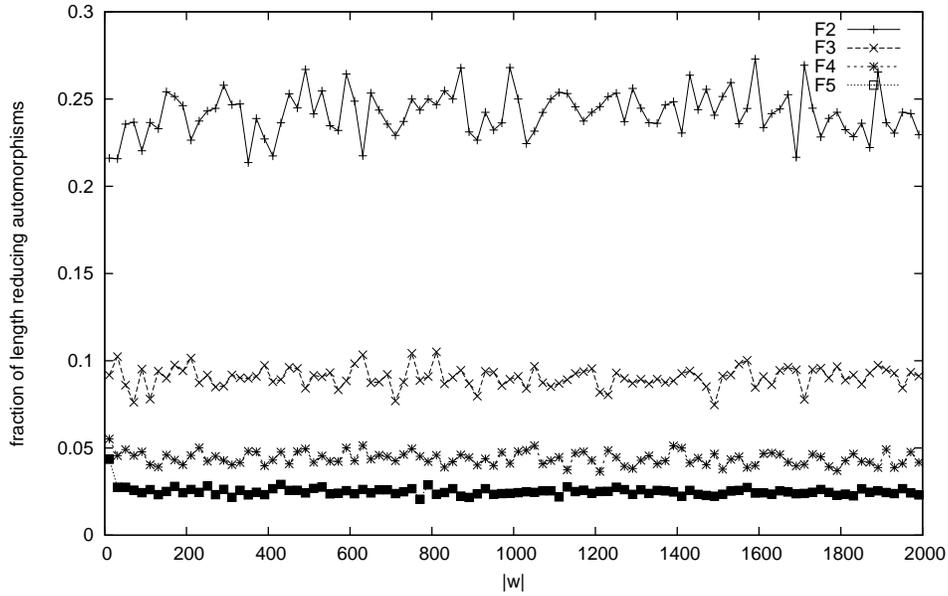}}
\caption{Estimates of $\frac{LR_n}{A_n}$ on inputs of various
lengths.} \label{fig:LR}
\end{figure}

The results of our experiments (for $k = 1$) indicate that the
average value of $|LR(U)|$ for a non-minimal $U$ of the total
length $l$ rapidly converges to a constant $LR_n$ when $l
\rightarrow \infty$.  In Table \ref{tab:LR}  and Figure
\ref{fig:LR}  we present  values of the $\frac{LR_n}{A_n}$ that
occur in our experiments for $k = 1$. This allows us to make the
following statement.
\begin{conclusion}
\label{con:DWA} The average number of length reducing Whitehead's
automorphisms for a given "generic" non-minimal  word $w \in F_n$
does not depend on the length of $|w|$, it depends only on the
rank $n$ of the free group $F_r$ (for sufficiently long words
$w$).
\end{conclusion}
A precise formulation of this statement is given  in Section
\ref{se:problems}.

%
\begin{table}[ht]
\begin{center}
\begin{tabular}{|l|c|c|c|c|}
\hline
 $|w|$                          & $F_2$ & $F_3$  &$F_4$ &$F_5$\\
\hline
 0..199  &  0.24 & 0.09 & 0.04 & 0.03 \\
\hline
 200..599  & 0.24 & 0.09 & 0.05 & 0.03\\
\hline
 600..999 & 0.24 & 0.09 & 0.04 & 0.02 \\
\hline
 1000..1299 & 0.25 & 0.09 & 0.04 & 0.02 \\
\hline
1400 ... 1800 & 0.24 & 0.09 & 0.04 & 0.02 \\
\hline
\end{tabular}
\end{center}
\caption{Estimates of $\frac{LR_n}{A_n}$ on inputs of various
lengths.} \label{tab:LR}
\end{table}
%

\subsubsection{ Algorithm DWB  }
\label{sub-sec:comp-DWB}

The obvious upper bound for the  time-complexity of  DWB
 is much higher, since one has to take into account all Whitehead automorphisms.
It is easy to see that there are
 $$B_n = 2n(2n-2)(2n-4) \ldots 2 = 2^n (n!) $$
 Whitehead automorphisms of the type 1).

 To run LS routine on $U$ it requires at most $d(A_n + B_n)$ runs of SR
 (which has complexity $c|U|$),  where $d$ is a constant which does not
 depend on $U$ and $n$. Now, to construct $Orb_{min}(U)$ it
 takes at most $|Orb_{min}(U)|$ runs of LS, hence one can bound
 the time complexity of DWA from above by
 $$d\cdot (A_n + B_n)\cdot c \cdot |U|\cdot |Orb_{min}(U)| .$$
This shows that DWB may be  very slow (in the worst-case) just because
there are too many Whitehead automorphisms in the rank $n$ for big $n$.
Moreover, the size of $Orb_{min}(U)$ can make the situation even worse.
Obviously,
 \begin{equation}
 \label{eq:orb-size}
 |Orb_{min}(U)| \leq 2n(2n-1)^{|U|-1},
\end{equation}
  hence a very rough estimates give the following upper bound for the
time-complexity of DWB:
$$d\cdot c \cdot (2n4^{(n-1)} - 2n + 2^n n!)\cdot |U| \cdot 2n(2n-1)^{|U|-1}.$$
One can try to improve on this upper bound through better
estimates of $ |Orb_{min}(U)|$. It has been shown in \cite{MS}
that for $k = 1$ and $n =2$ the number $ |Orb_{min}(U)|$ is
bounded from above by a polynomial in $|U_{min}|$. It was also
conjectured in \cite{MS} that this result holds for arbitrary $n
\geq 2$,  and for $n = 2$ the upper  bound is the following:

 $$|Orb_{min}(U)| \leq 8|U_{min}|^2 +40|U_{min}|.$$
Recently, B.Khan proved in \cite{Khan} that the bound above holds,
indeed. Still, independently of the size of the set
$Orb_{min}(U)$, the number $B_n$ of elementary Whitehead
automorphisms in rank $n$ makes DWB impractical for sufficiently
big $n$.

The net outcome of the discussion above is that the  algorithms
DWA and DWB  are intractable for "big"  ranks, even though for a
fixed rank $n$ DWA is quadratic in $|U|$  and DWB could  be
polynomial in $|U|$ (if Conjecture \ref{conj:1} from Section
\ref{se:problems} holds).

\subsection{General Length Reduction Problem.} \label{sec:LRP}

Observe that the main part of DWA is the elementary length reduction
routine ELR, which for a given tuple $U \in F^k$ finds a {\it
Whitehead} automorphism $\varphi \in \Omega(F)$ such that
\begin{equation} \label{eq:1}
|U\varphi| < |U|
\end{equation}
An arbitrary automorphism $\varphi \in Aut(F)$ is called {\it
length-reducing} for $U$ if it satisfies the condition
(\ref{eq:1}) above.

Obviously,  to solve Problem A it  suffices to find an arbitrary (not
necessary Whitehead) length-reducing automorphism for a non-minimal
tuple $U$. We have seen in Section \ref{sub-sec:est} that the
time-complexity of the standard Whitehead algorithm for Problem A
depends mostly on the cardinality of the set $\Omega_n$ which is huge
for big $n$. One of the key ideas on improving the efficiency of
Whitehead algorithms is to replace $\Omega_n$ by  another smaller set
of automorphisms of $F$ or to use a different strategy to find
length-reducing automorphisms.  To this end we formulate the following

\medskip {\bf Length-Reduction Problem (LRP).}
{\it For a non-minimal  tuple $U \in F^k$ find a length-reducing
automorphism.}

\medskip

Theorem \ref{T1} gives one solution to LRP - the algorithm DWA. In
Section \ref{sec:GA} we describe a genetic algorithm which, we believe,
solves LRP much more efficiently on average then DWA.

\section{Description of the genetic  algorithm}

\label{sec:GA}

In this section we describe Genetic Whitehead Algorithm (GWA) for
solving Whitehead's  Problem A.

Genetic algorithms are stochastic search algorithms driven by a
heuristic, which is  represented by an evaluation function, and
special random operators: crossover, mutation and selection.

Let $\mathcal{S}$ be a search space. We are looking for an element
in $\mathcal{S}$ which is a solution to a given problem. A tuple
$P \in \mathcal{S}^r$  ($r$ is a fixed positive integer) is called
a {\it population}   and components of $P$ are called {\it
members} of the population.  The initial population $P_0$ is
 chosen randomly.
 On each {\it iteration}  $i = 1, 2, \ldots$ Genetic Algorithm produces a new
 population $P_i$ by means of random operators.  The
goal is to produce a population which contains a solution to the
problem.   One iteration of Genetic Algorithm  simulates natural
evolution. A so-called  {\it fitness function} $Fit : \mathcal{S}
\rightarrow \mathbb{R}_+$ implicitly directs this evolution:
members of the current population $P_i$ with higher fitness value
have more impact on generating the next population $P_{i+1}$. The
function $Fit$
 measures on how close is a given member $m$ to a
solution. To halt the algorithm   one has to provide in advance  a
 {\it termination condition} and check whether it holds or not on each
 iteration.   The basic
structure of the standard Genetic Algorithm is given in Figure
\ref{fig:GA}.

\begin{figure}[ht]
\centerline{\rule{12cm}{0.5pt}}
\medskip
\emph {{\bf procedure} Genetic Algorithm \\
Initialize current  population $P \in \mathcal{S}^r$; \\
Compute fitness values $Fit(m)$, $ \forall m \in P$; \\
{ \bf WHILE NOT }  the termination condition satisfied {\bf DO }
\begin{list}{}{}
\item If we assume that greater values of function $Fit$ correspond to
the better solutions, then the probability $Pr(m)$ of the
member $m \in P$ to be selected
\[ Pr(m) = \frac{Fit(m)}{ \sum_{m_i \in P} Fit(m_i)},\]
\item Create new members by applying crossover and/or mutation to
the selected members; \item Generate a new population by replacing
members of the current population by the new ones;
 \item Recompute fitness values;
\end{list}
{\bf END WHILE LOOP } }

\medskip
\centerline{\rule{12cm}{0.5pt}}
 \label{fig:GA} \caption{  Structure of the
standard Genetic Algorithm}
\end{figure}

The choice of random operators and evaluating functions is crucial
here. This requires some problem specific knowledge and a good
deal of intuition. Below we give detailed description of the major
components of the genetic algorithm GWA for solving Problem A.

\subsection{Solutions and members of the  population}
\label{sec:Sol}

Solutions to the  Problem A are finite sequences of Whitehead
automorphisms which carry a given tuple $U \in F^k$ to a minimal tuple
$U_{min}$. As we have mentioned above one may use  only automorphisms
of the type 2) for this problem. Moreover, not all automorphisms of the
type 2) are needed as well (recall that a big number of such
automorphisms is the main obstacle for the the standard Whitehead
algorithm DWA).  What are optimal sets of  automorphisms is an
interesting problem which we are going to address in \cite{HMM}, but
our preliminary experiments show that the following set gives the best
results up to date.

Let  $X = \{x_1,...,x_n\}$  and $F = F_n(X)$.  Denote by $T =T _n$ the
following set of Whitehead automorphisms:\\
$(W1) \;x_{i} \rightarrow x_{i}^{-1},\: x_{l} \rightarrow x_{l},$ \\
$(W2) \;x_{i} \rightarrow x_{j}^{\pm 1}x_{i},\: x_{l} \rightarrow x_{l},$ \\
$(W3) \;x_{i} \rightarrow x_{i}x_{j}^{\pm 1},\: x_{l} \rightarrow x_{l},$\\
$(W4) \;x_{i} \rightarrow x_{j}^{-1}x_{i}x_{j},\: x_{l} \rightarrow
x_{l},$\\
 where $i \neq j$ and  $i \neq l$.

We call $T$ the \textit{ restricted set} of Whitehead transformations.
It follows from \cite{Nielsen} that  $T$ generates $Aut(F)$. Hence any
solution to Problem A can be represented by a finite sequence of
transformations from $T$.
 Notice that $T$ has much fewer elements than $\Omega _n$:
\[|T| = 5n^2 - 4n.\]

We define the search space $\mathcal{S}$  as the set of all finite
sequences $\mu  = <t_1, \ldots, t_s>$ of transformations from $T$. For
such $m$ and a tuple $U \in F^k$ we define $U\mu = Ut_1 \ldots t_s$.

At the beginning the algorithm  generates an initial population by
randomly selecting members. How to choose the size of the initial (and
all other) population is a non-trivial matter. It is clear that bigger
the size - larger the search space which is explored in one generation.
But the trade off is that we may be spending too much time evaluating
fitness value of members of the population. We do not know the optimal
size of the population, but populations with 50 members seem to give
satisfactory results.

\subsection{Evaluation methods}
\label{sec:Fit}

\emph{Fitness function} $Fit$  provides a  mechanism to
  assess  members of a given population $P$.

Recall that the aim of GWA is to find a sequence of transformations
$\mu = ( t_1, \ldots, t_s ), \ t_i \in T,$ such that $$U\mu = U_{min}$$
for a given input $U  \in F^k$. So members $\mu$ of a given population
$P$ with  smaller total length $|U\mu|$ are closer  to a solution,
i.e., "fitter", than the other members.  Therefore we define the
fitness function $Fit$ as
\[Fit(\mu) = \max_{\lambda \in P}\{|U\lambda|\} - |U\mu|.\]
Observe, that members with higher fitness values are closer to a
solution $U_{min}$ with respect to the metric on the graph
$\Gamma(F,k,T)$. In fact, we have two different implementations of
the evaluation criterion: the one as above, and another one in
which a word is considered as a {\em cyclic} word, so we evaluate
fitness values of  cyclic permutations of $U\lambda$.

\subsection{Termination condition}
\label{sec:AG}

 {\it Termination condition} is a tool to check whether a
given population contains a solution to the problem or not.

In the  case of Whitehead method there are several ways to define a
termination condition.

T1) Once a new population $P_n$ has been defined and all members of it
have been evaluated one may check whether or not  $P_n$ contains a
solution to Problem A.  To this end one can run "Elementary Length
Reduction" routine on $U\mu*$ for each fittest member $\mu^* \in P_n$
until $U_{min}$ is found.  Theoretically, it is a good termination
condition, but, as we have mentioned already,  to run ELR might be very
costly.

T2)  If  for a given tuple $U$ we know in advance   the length of a
minimal tuple $|U_{min}|$ ( for example, when  $U$ is a part  of a
basis of $F$), then we define another (fast) termination condition as
$|U\mu^*| = |U_{min}|$ for some fittest member $\mu^* \in P_n$.

T3) Suppose now that  we do not know $|U_{min}|$ in advance, but we
know the expected number of populations, say $E = E(U)$, (or some
estimates for it) which is required for the genetic algorithm GWA to
find $U_{min}$ when starting on a tuple $U$. In this case we can use
the following strategy: if the algorithm keeps working without
improving on the fitness value $Fit(\mu^*)$ of the fittest members
$\mu^*$  for long enough,  say for the last $pE$ generations (where $p
\geq 1$ is a fixed constant), then it halts and gives $U\mu^*$ for some
fittest $\mu^*$ as an outcome.

If the number $E = E(U)$ is sufficiently small this termination
condition could be efficient enough. Below, we will describe some
techniques and numerical results on how one can   estimate the number
$E(U).$ Of course, in this case there is no guarantee that the tuple
$U\mu^*$ is indeed minimal.  We refer to  such termination conditions
as to {\it heuristic } ones, while the condition T1 is {\it
deterministic}.

T4) One can combine conditions T3 and T1 in the following way. The
algorithm uses the heuristic termination condition T3 and then checks
(using T1) whether or not  the output $U\mu^*$  is indeed minimal. It
is less costly then T1 (since we do not apply T1 at every generation)
and it is more costly then T3.

\subsection{Stochastic operators}
There are five basic random operators that where used in the
algorithm.

\subsubsection{One point crossover}

 Let  $\mu_{1} = <t_1,...,t_e>$ and $\mu_{2}= <s_1,...,s_l>$ be two
members of a population $P_n$ which are choosen with respect to some
selection method.  Given two random numbers $ 0 < p< e$ and $0 < q < l$
the algorithm  constructs two offsprings $o_1$ and $o_2$ by {\it
recombination} as follows:
\begin{eqnarray*}
o_{1}=<t_1,...,t_{p-1},s_q,...,s_l>, \ \ \
o_{2}=<s_1,...,s_{q-1},t_p,...,t_e>.
\end{eqnarray*}

\subsubsection{Mutations}

The other four operators $M_{att}, M_{ins}, M_{del}, M_{rep}$ act on a
single member of a population and are usually called { \em mutations }.
They attach, insert, delete, or replace some transformation in a
member. Namely, let $\mu = <t_{1},...,t_{l}>$ be a member of a
population. Then:

\begin{enumerate}
\item [$M_{att}$] \textit{attaches} a random transformation $s\in T$
\[M_{att} :  <t_{1},...,t_{l}> \; \rightarrow \; <t_{1},...,t_{l},s>;\]
\item [$M_{ins}$] \textit{inserts} a random transformation $s\in T$ into a
randomly chosen position $i$
\[M_{ins} : <t_{1},...,t_{l}> \; \rightarrow \; <t_{1},...,t_{i-1},s,t_{i},..,t_{l}>;\]
\item  [$M_{del}$] \textit{deletes} the transformation in  a randomly chosen
position $i$
\[M_{del} : <t_{1},...,t_{l}> \; \rightarrow \; <t_{1},...,t_{i-1},t_{i+1},..,t_{l}>;\]
\item [$M_{rep}$] \textit{replaces} the randomly chosen $t_i$ by
a randomly chosen $s \in T$
\[M_{rep} : <t_{1},...,t_{l}> \; \rightarrow \; <t_{1},...,t_{i-1},s,t_{i+1},..,t_{l}>.\]
\end{enumerate}
Operator $M_{att}$ is a special case of $M_{ins}$, but it is convenient
to have it as separate operator (see remarks in the Section
\ref{sub-sec:3.5.1}).

\subsubsection{Replacement}
\label{sub:selection}

In this section we discuss  a protocol to construct members of the next
population $P_{new}$ from the current population $P$.

First, we select randomly two members $\mu, \lambda$ from $P$. The
probability to choose a member from $P$ is equal to
\[ Pr(m) = \frac{Fit(m)}{ \sum_{m_i \in P} Fit(m_i)}.\]
With small probability (0.10 - 0.15) we add both $\mu$ and $\lambda$ to
an intermediate population $P_{new}^{\prime}$. Otherwise, we apply the
crossover operator to $\mu$ and $\lambda$ and add the offsprings to
$P_{new}^{\prime}$. We repeat this step until we get the required
number of members in $P_{new}^{\prime}$ (in our case 50).

Secondly, to every member $m \in P^{\prime}_{new}$ we apply a random
mutation $M$ with probability 0.85 and add the altered member to the
new population $P_{new}$. The choice of $M$ is governed by the
corresponding probabilities $p_M$. Otherwise (with probability 0.15) we
add the member $m$ to $P_{new}$ unchanged.  We refer to Section
\ref{sub-sec:3.5.1} for a detailed discussion of our choice of the
probabilities $p_M$.

In addition the solution with the highest fitness value among all
previously occurred solutions is always added to the new population
(replacing a weakest one). This implies that if we denote by $\mu_n$
one of the fittest members of a population $P_n$ then
$$|U\mu_0| \geq |U\mu_1| \geq \ldots $$


\subsection{Some important features of the algorithm}
\label{sec:Rem}
\subsubsection{Precise solutions and local search}
\label{sub-sec:3.5.1}

It has been shown that different heuristics and randomized methods can
be combined together, often resulting in more efficient \emph{hybrid}
algorithms. Genetic algorithms are  good in covering large areas of the
search space. However, they may fail when a more thorough trace of a
local neighborhood is required. In case of symbolic computations this
becomes an important issue since we are looking for an exact solution,
not an approximate one.  Even if the current best member of a
population  is one step away from the optimum it might take some time
for the standard genetic algorithm to find it. In our case, experiments
show that the standard genetic algorithms can quickly reach the
neighborhood of the optimum, but it  may be stuck being unable to hit
the right solution. To avoid that one could add a variation of the
local search procedures to the standard genetic algorithm.

In GWA some kind of gradient descent procedure was implicitly
introduced via   mutation operators. Observe,  that in general, if $M
\neq M_{att}$ then  for a given member $\mu$ the tuple $UM(\mu)$ lies
far apart from $U\mu$ in the graph $\Gamma(F,k,T)$.  However, the
mutation $M_{att}$  always gives a tuple $U M_{att}(\mu)$ at distance 1
from $U \mu$ in the graph $\Gamma(F,k,T)$. Therefore, the greater
chance to apply $M_{att}$,  the more neighbors of $U\mu$ we can
explore. It was shown experimentally that GWA performs much better when
$M_{att}$ has a greater chance to occur.  We used $p_{M_{att}} = 0.7$,
and $p_M = 0.1$ for $M \neq M_{att}.$

\subsubsection{Substitution Method}
\label{sub:substitution}

One of the major concerns when dealing with a search problem is that
the algorithm may fall into a local minimum. Fortunately, Theorem
\ref{T1} shows that every local minimum of the fitness function $Fit$
is, in fact,  a global one. This allows one to introduce another
operator, which we call \textit{Substitution}, and which is used  to
speed up the convergence of the algorithm.

Suppose that the algorithm  found  a member $\mu_n \in P_n$ which is
fitter than all the members of the previous population $P_{n-1}$ (a
genetic variation of ELR routine). Then we want our algorithm to focus
more on the tuple $U\mu$ rather then  to spread its own resources for
useless search elsewhere. To this end, we stop the algorithm and
restart it replacing the initial tuple $U$ with the tuple $U\mu$ (of
course, memorizing the sequence $\mu$). That is a genetic variation of
the Whitehead's gradient descent (see Section \ref{sub-sec:WA}).
 This simple method has tremendously improved the
performance of the algorithm. In a sense,  this  substitution turns GWA
into an algorithm which  solves a sequence of  Length Reduction
Problems.


\section{Experiments and Results}
\label{sec:results}

 Let $F = F_r(X)$ be a free group of rank $r$ with basis
$X$. For simplicity we describe here only experiments with
Whitehead algorithms on inputs from $F$ (not arbitrary $k$-tuples
from $F^k$). Moreover, in the present paper we focus only on the
time-complexity of Problem A, leaving discussion on Problem B for
the future. In fact, we discuss mostly the length reduction
problem LRP, as a more fundamental problem. In our experiments we
choose ranks $r = 2,5,10,15,20$.
   Before we going into details it is worthwhile to discuss a few
 basic problems on statistical analysis of experiments with
 infinite groups.

\subsection{Experimenting with infinite groups.} \label{sec:Exp}

In this section we discuss briefly several general problems
arising in experiments with infinite groups.

Let ${\cal A}$ be an algorithm for computing with elements from a
free group $F = F_r(X).$ Suppose that  the set of all possible
inputs for ${\cal A}$ is an infinite subset $S \subset F$.
Statistical analysis of experiments with  ${\cal A}$  involves
three basic parts:
 \begin{itemize}
 \item creating a finite set of test inputs
$S_{test} \subset S$,
 \item running ${\cal A}$  on inputs from $S_{test}$
and collecting outputs,
 \item statistical analysis of the  resulting
data. \end{itemize}
 The following is the main concern when creating  $S_{test}$.

\medskip
{\bf Random Generation of the  test data: } {\it How one can
generate pseudo-randomly a finite subset $S_{test} \subset S$
which represents adequately the whole set  $S$?}

\medskip
\noindent The notion of a random element in $F$, or in $S$,
depends on a chosen  measure on $F$. Since $F$ is infinite,
elements in $F$ are not  uniformly distributed. The problem cannot
be solved just by replacing $F$ with a finite ball $B_n$,  of all
elements in $F$ of length at most $n$, for a big number $n$.
Indeed, firstly,  the ball $B_n$ is too big for any practical
computations; secondly,  from group-theoretic view-point elements
in $B_n$ usually are not uniformly distributed.   We refer to
\cite{BMS} and \cite{BMR} for a thorough discussion of this
matter.

The main problem  when collecting  results of the runs of the
algorithm ${\cal A}$ on inputs from $S_{test}$   is pure
practical: our resources in time and computer power are limited,
so the set $S_{test}$ has to be as small as possible, though still
representative.

\medskip
{\bf Minimizing the cost: } {\it How to make the set $S_{test}$ as
small as possible, but still representative? }

\medskip
\noindent Below we used the following  technique  to ensure
representativeness of  $S_{test}$.    Assume we have already a
procedure to generate pseudo-random elements in $S$. Let
$\chi(S_{test})$ be some computable numerical characteristic of
the set $S_{test}$, which represents a "feature" that we are going
to test. Fix  a small real number $\varepsilon
> 0$. We start creating $S_{test}$ by generating an
initial subset $S_0 \subset S$ which we can easily handle within
our recourses. Now we enlarge the set $S_{0}$ to a new set $S_1$
by pseudo-randomly adding reasonably many of new elements from
$S$, and check whether the equality
$$|\chi(S_0) -
\chi(S_1)| \leq \varepsilon$$ holds or not. We repeat this
procedure until the equality holds for $N$ consecutive steps $S_i,
S_{i+1}, \ldots, S_{i+N}$, where $N$ is a fixed preassign number.
In this event we stop and take  $S_{test} = S_i$.

Statistical analysis of the experiments depends on the  features
that are going to be tested (average running time of the
algorithm, expected frequencies of outputs of a given type, etc.).
For example, estimations of the running time of the algorithm
${\cal A}$  depends on how we measure  "complexity" or "size" of
the inputs $s \in S$. For example, it turned out that the running
time of the Whitehead algorithm GWA does not depend essentially on
the length of an input word $s$, so it would be meaningless  to
measure the time complexity of $DWA$  in terms of the  length of
$s$ (as is customary in computer science). So the following
problem is crucial here.

\medskip
{\bf Finding adequate complexity functions: } {\it Find a
complexity function on $S$ which is compatible with the algorithm
${\cal A}$. }

\medskip
\noindent

Below we suggest some particular ways to approach all these
problems in the case of the Whitehead's algorithms.

\subsection{Random elements in $F$ and Whitehead algorithms}
\label{sub-sec:random}

 It seems that the most obvious choice for the set $S_{test}$ to test
 performance of various Whitehead algorithms would be a
  finite set $S_F$ of randomly chosen elements from $F$. It turned out,
   that this choice is not good at all since with a high
  probability a random element in $F$ is already  minimal. Nevertheless,
  the set  $S_F$ plays an important part in the sequel as a base for
other  constructions.

 A random element $w$  in $F = F_r(X)$ can be produced as the result of a
no-return simple random walk on the Cayley graph of $F$ with
respect to the set of generators $X$ (see \cite{BMR} for details).
In practice this amounts to a pseudo-random choice of a  number
$l$ (the length of $w$),  and  a pseudo-random  sequence $y_1,
\ldots, y_l$ of elements $y_i \in X^{\pm 1}$ such that $y_i \neq
y_{i+1}^{-1}$, where $y_1$ is chosen randomly from $X^{\pm 1}$
with probability $1/2r$, and all others are chosen randomly with
probability $1/(2r-1)$.  It is convenient to structure the set
$S_{F}$ as follows:
$$ S_F =
\bigcup_{l = 1}^L S_{F,l}, \ \ \  S_{F,l} = \bigcup_{i = 1}^K
w_{i,l}$$ where $w_{i,l}$ is a random word of length $l$ and  $L,
K$ are parameters.

\begin{figure}[ht]
\centerline{  \includegraphics{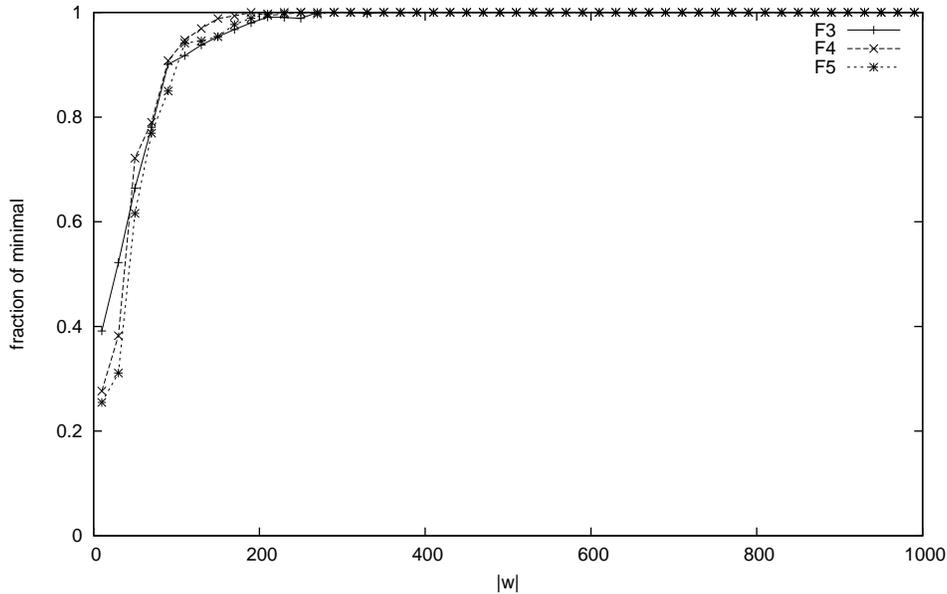}} \caption{
Fractions of Whitehead-minimal elements in a free group $F_r$, $r=
3,4,5$.} \label{fig:SF}
\end{figure}

 To find all minimal elements in $S_F$  we run
 the standard deterministic Whitehead algorithm DWA on every $s\in
 S_F$. Since DWA is very slow for big ranks we experimented with
 free groups $F = F_r$ for $r = 3,4,5$.
 In Figure \ref{fig:SF} we present the fractions of minimal
elements among all elements of a given length in $S_F$.

%
%

This experimental data leads to the following statement.
\begin{conclusion}
Almost all elements in $F_r, r \geq 2$ are Whitehead minimal.
\end{conclusion}
We refer to Section \ref{se:problems} for a rigorous formulation
of the corresponding mathematical  statement.

 The running time  $T_{DWA}(w)$ of  the standard Whitehead algorithm DWA on
  a minimal input $w$ is very easy to estimate. Indeed, in this case DWA
  applies the substitution routine SR for every Whitehead automorphism of
the second type.
   Since there are $A_r$ such automorphisms (see Section \ref{sub-sec:WA}), then

$$ A_r \leq T_{DWA}(w) \leq c \cdot A_r|w|.$$

The time spent by the  genetic  algorithm GWA on a random input
$w$ depends solely on the  build-in  termination condition: if it
is heuristic (see Section \ref{sec:AG}),   then  GWA stops after
$pE(w)$ iterations, where $E(w)$
 is the expected running time for GWA on the  input $w$; if it is deterministic
 then again it takes  $A_r$ steps for GWA to halt.
This shows that the set $S_F$ does not really test how GWA works,
instead, it  tests only the termination conditions.

We summarize the discussion above in the following statement.

\begin{conclusion}
The time-complexity of Whitehead algorithms DWA and GWA on generic
inputs from  $S_F$ is easy to estimate. The set $S_F$ does not provide
any means to compare algorithms DWA and GWA.
\end{conclusion}
It follows that one has to test Whitehead algorithms on inputs $w \in
F$  which are {\it non-minimal}.

\subsection{Complexity of Length Reduction Problem}
\label{sec:WC1}

In this section we test our genetic algorithm GWA  on the length
reduction problem LRP, which is the main component  of the
Whitehead's Method.

To this end we  generate a finite set $S_{NMin}(r)$ of non-minimal
elements in a free group  $F_r$, for $r = 2,5,10,15,20,$ by
applying random Whitehead automorphisms to elements form $S_F$.
More precisely, put
\[ S_{NMin}(r) = \bigcup_{l} \bigcup_{1 \leq i \leq K} w_{i,l} \varphi_i ,\]
where  $\varphi_i$ is a randomly chosen Whitehead automorphism of
type 2), $w_{i,l} \in S_F$ with $|w_{i,l}| < |w_{i,l} \varphi_i|$.
Since almost all elements from $S_F$ are minimal it is easy to
generate a set like $S_{NMin}(r)$.  Notice that elements in
$S_{NMin}(r)$ are not randomly chosen  non-minimal elements from
$F$, they are non-minimal elements at distance 1 from minimal
ones. We will have to say more about this in the next section.

The results of our experiments indicate that the average time
required for GWA to find a length reducing Whitehead automorphism
for a given non-minimal element $w \in S_{NMin}(r)$ does not
depend significantly on the length of the word $w$.

Let $T_{gen}(w)$ be the number of iterations required for GWA to
find a length-reducing automorphism for a given $w \in F$ during a
particular run of GWA on the input $w$. We compute the average
value of $T_{gen}(w)$ on inputs $w \in S_{NMin}(r)$ of a given
"size". If the length of a word $w$ is taken as its size then  we
obtain the following time complexity function with respect to the
test data $S_{NMin}(r)$:
$$T_r(m) = \frac{1}{|S_m|}\sum_{w \in S_m}T_{gen}(w)$$
where $S_m = \{w \in S_{NMin}(r) \mid |w| = m\}$.

\begin{figure}[ht]
\centerline{  \includegraphics{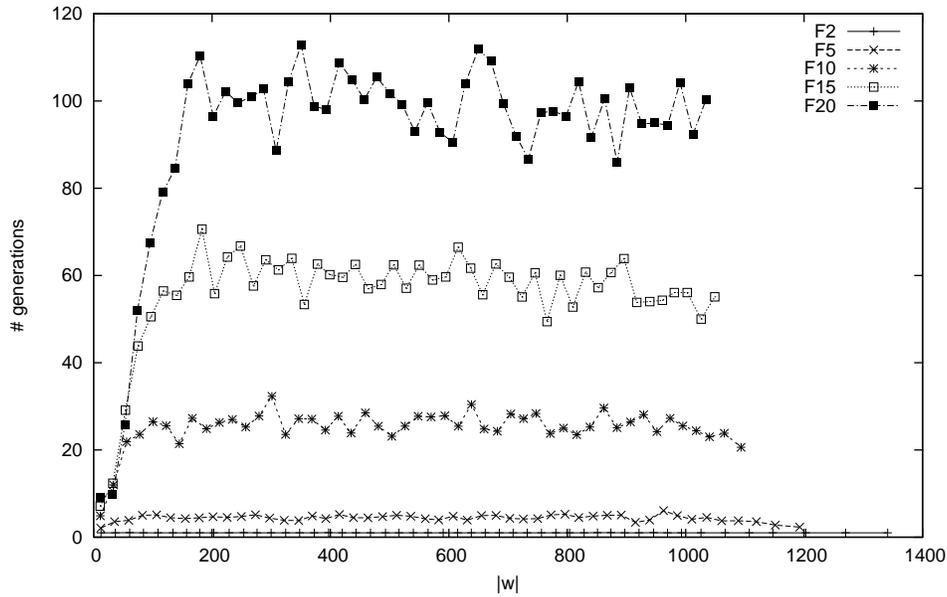}} \caption{Values of
$T, S = S_{1}$.} \label{fig:WGA_non_min}
\end{figure}

Values of $T_r(m)$ are presented  in Figure \ref{fig:WGA_non_min}
for free groups $F_r$ with $r = 2,3,5,10,15,20.$

 We can see from
the graphs that the function $T_r$ grows for small values of $|w|$
and then stabilizes at some constant value $T_r^\ast$. This shows
that $T_r$ does not depend  on the word's length and depends only
on the rank $r$   (for long enough words $w$).

 In  Table \ref{tab:correl} we give correlation coefficients
between $T_r$ and $|w|$  for $r = 2,5,10,15,20$, which are
sufficiently small.


%
\begin{table}[ht]
\begin{center}
\begin{tabular}{|l|c|c|c|c|c|c|c|c|c|}
\hline
            &  $F_2$ & $F_5$ & $F_{10}$ &  $F_{15}$ & $F_{20}$\\
\hline
all words   & -0.012 & -0.016  & 0.015  &  0.03  & 0.072\\
\hline
$|w| > 100$ & -0.011  & -0.03 & -0.019  & -0.025  & -0.005\\
\hline
\end{tabular}
\end{center}
\caption{Correlation between $|w|$ and $T_r$.}
 \label{tab:correl}
\end{table}
%

 We summarize the discussion above in the following
 statements.
 \begin{conclusion}
 \label{con:2}
 The number of iterations required for GWA to find a length reducing automorphism
for  a given non-minimal input $w$ does not
 depend on the length of $|w|$,  it depends only
 on the rank $r$  (for long enough input words).
 \end{conclusion}
Recall that a similar phenomena was observed for the deterministic
Whitehead's algorithm in Conclusion \ref{con:DWA}.
\begin{conclusion}
 \label{con:3}
One has to replace the length size function by a more sensitive "size"
function  when measuring the time-complexity of the Length Reduction
Problem.
\end{conclusion}

\begin{conclusion}
 \label{con:4}
For each free group $F_r$ the time-complexity function $T_r$ is
bounded from above by  some constant value $T_r^\ast$.
\end{conclusion}
 We can try to estimate the value $T_r^\ast$ as  the {\it expected number of
generations}
 $$E(r) = \frac{1}{|S_{NMin}(r)|}\sum_{w \in S_{NMin}(r)}T_{gen}(w).$$
  required for GWA to find a length-reducing automorphism for
generic non-minimal elements from $F_r$. Notice, that we use
$E(r)$ in the  heuristic termination condition TC3 (see Section
\ref{sec:AG}) for the algorithm GWA.

Of course, the conclusions above are not mathematical theorems,
they are just empirical phenomena that can be seen from our
experiments based on the test  set $S_{NMin}(r)$. It is important
to make sure that the set $S_{NMin}(r)$ is sufficiently
representative.

To this end, we made sure, firstly, that the distributions of
lengths of words from the set $S_{NMin}(r)$ are similar for
different ranks (using the variable  $l$). Secondly, our choice of
the parameter $K$ in the construction of $S_{NMin}(r)$  ensures
representativeness of the test data with respect to the
characteristic $E(r)$. Namely, we select $K$ such that for larger
values  $K^\prime > K$ the corresponding value $E_{K^\prime}(r)$
does not differ significantly from $E_K(r)$ (here $E_K(r)$ is the
value corresponding to the data set $S_{NMin}(r)$ with the
parameter $K$).

Values of $E(r)$ for  different $K$ and $r$   are given in Table
\ref{tab:S_non_min_repres}.

%
\begin{table}[ht]
\begin{center}
\begin{tabular}{|l|c|c|c|c|c|}
\hline
 $K$    &    $E(2)$ & $E(5)$ & $E(10)$ & $E(15)$ & $E(20)$  \\
\hline
100 &  1.007 & 2.43 & 6.55 & 11.48 & 16.98  \\
\hline
200   & 1.009 & 2.42 & 6.44 &  11.47 & 17.17 \\
\hline
300   &  1.008 & 2.42 & 6.43 & 11.39 & 17.3\\
\hline
400 &   1.007 & 2.39 & 6.43 & 11.40 &  17.38 \\
\hline
500 &   1.007 & 2.44 & 6.43 & 11.39 & 17.4 \\
\hline
\end{tabular}
\end{center}
\caption{$E_K(r)$ for  different values of $K$ and $r$. }
 \label{tab:S_non_min_repres}
\end{table}

\subsection{Complexity functions}
\label{subsec:complexity}

In this section we discuss possible complexity, or size, functions
suitable to estimate the time-complexity of different variations
of Whitehead algorithms.  Below we suggest a new complexity
function based on the distance in the Whitehead graph.

Let $F = F_r$, $Y  \subset Aut(F)$  a set of generators of the
group $Aut(F)$,   $\Gamma(F,Y) = \Gamma(F,1,Y)$ the Whitehead
graph on $F$  relative to $Y$ (see Section \ref{sub-sec:WTh}). For
a word $w \in F$ we define $WC_Y(w)$ as a minimal number of
automorphisms from $Y^{\pm 1}$ required to reduce $w$ to a minimal
one $w_{min}$. Notice that $WC_Y(w)$ is the length of a geodesic
path in $\Gamma(F,Y)$ from $w$ to some $w_{min}$. If $Y$ is the
set of all Whitehead automorphism $\Omega_r$ then we call
$WC_Y(w)$ the {\it Whitehead's complexity} of $w$ and denote it by
$WC(w)$. Similarly, one can introduce the Nielsen's complexity of
$w$, $T$-complexity, etc. In this context minimal elements have
zero Whitehead complexity.

\medskip
{\bf Claim } {\em The Whitehead's complexity function $WC(w)$ is
an adequate complexity function to measure performance of various
modifications of the Whitehead's algorithm.}

\medskip\noindent
Indeed, let ${\cal K}$ be a class of Whitehead's-type algorithms
which use an arbitrary  generating set $Y \subset \Omega_r$ of
Whitehead automorphisms to find a minimal word $w_{min}$ for an
input word $w$. The best possible algorithm of this type is
 the non-deterministic Whitehead
algorithm NDWA with an oracle that at each step $i$ gives a length
reducing automorphism $t_i \in Y$ such that $|wt_1\ldots t_i| <
|wt_1\ldots t_{i-1}|$.  Clearly, it takes $WC_Y(w)$  steps for
NDWA to produce  $w_{min}$. Thus, measuring efficiency of an
algorithm $A \in {\cal K}$ in terms of $CW_Y$ gives us a
comparison of performance of ${\cal A}$ to the performance of the
best possible  algorithm in the class.

\begin{remark}
Notice that the set $S_{NMin}(r)$ is a pseudo-random sampling of
elements $w \in F_r$ with $WC(w) = 1$. This explains the behavior
of the function $T_r$ in  Figure \ref{fig:WGA_non_min}. The number
of iterations required for GWA to find a length reducing
automorphism depends on Whitehead complexity not on the lengths of
the words.
\end{remark}

 Of course, $WC$ complexity is mostly a theoretical tool, since,
 in general, it is harder  to compute $WC(w)$ then to find $w_{min}$.
It follows from the Whitehead's fundamental  theorem that $WC(w)
\leq |w|$ for every $w \in F$. In Table \ref{tab:CW-W} we
collect some experimental results on relation between $WC(w)$ and
$|w|$.

%
\begin{table}[ht]
\begin{center}
\begin{tabular}{|l|c|c|c|c|c|c|}
\hline
                           & $F_2$  &$F_5$ & $F_{10}$ & $F_{15}$ & $F_{20}$\\
\hline
 $|w t| / |w|, \; t \in \Omega$ & 1.04 & 1.20  & 1.26 &  1.28 & 1.29 \\
\hline
 $|w t| / |w|, \; t \in T$ & 1.06 &  1.15 &  1.10 & 1.07  & 1.06 \\
\hline
\end{tabular}
\end{center}
\caption{$WC(w)$ vs $|w|$.}
\label{tab:CW-W}
\end{table}
%

This leads to the following
\begin{conclusion}\label{con:cr}
Let $W_m = \{w \in F_r \mid WC(w) = m\}$. Then there exists a
constant $c_r$ such that $$|w| \geq c_r^m$$ for the "most"
elements in $W_m$.
\end{conclusion}

For the stochastic algorithm GWA one can define an average time
complexity function $T_{r,Y}(m)$ with respect to the test data
$S_{NMin}(r)$ and the "size" function $WC_Y$ as follows:
$$T_{r,Y}(m) = \frac{1}{|S_m|}\sum_{w \in S_m}T_{gen}(w)$$
where $S_m = \{w \in S_{NMin} \mid WC_Y(w) = m\}$.

\begin{conjecture} \label{con:WCcomp}
The average number of iterations required for GWA to find
$w_{min}$ on an input $w \in F$ depends only on $WC(w)$  and the
rank of the group $F$.
\end{conjecture}
We discuss some experiments made to  verify Conjecture
\ref{con:WCcomp} in Section \ref{sec:prim}.

\subsection{Experiments with primitive elements}
\label{sec:prim}

In this section we discuss results of experiments with primitive
elements. Recall that elements from the orbit $Orb(x_i)$, where
$x_i \in X$, are called {\it primitive} in $F(X)$. Experimenting
with primitive elements has several important advantages:
 \begin{itemize}
 \item in general,  primitive elements $w$ require  long
 chains of Whitehead automorphisms (relative to $|w|$) to get to
$w_{min}$,
 \item one can easily generate pseudo-random primitive elements,
 \item the genetic algorithm GWA has a perfect termination condition
 $|w_{min}| = 1$  for
 primitive elements $w$.
 \end{itemize}
 Thus, primitive elements provide an optimal test data to compare various
 modifications of Whitehead algorithm and to verify (experimentally)
 the conjectures and conclusions stated in the previous sections.

We generate primitive  elements in the form $x \varphi$,  where
$x$ is a random  element from $X$ and $\varphi$ is a random
automorphism of $F$ given by a freely reduced product $\varphi =
t_1 \ldots t_l$  of $l$ randomly and  uniformly chosen
automorphisms from $T$ with $t_i \neq t_{i+1}^{-1}$ (see the
comments for $S_F$). The number $l =l(\varphi)$ is called the
 {\it length } of $\varphi$.

 In general, a  random automorphism $\varphi$ with respect to a
 fixed finite set $T$  of
generators of the group $Aut(F)$  can be generated as the result
of a no-return simple random walk on the Cayley graph
$\Gamma(Aut(F),T)$ of $Aut(F)$ with respect to the set of
generators $T$. Unfortunately, the structure of $\Gamma(Aut(F),T)$
is very complex, and it is hard  to simulate such a random walk
effectively.

Again, for each free group $F_r$ ($r = 2,5,10,15,20$), we construct a
set $S_P(r)$ of test primitive elements as follows:
\[ S_P(r) = \bigcup_{l=1}^{L} \bigcup_{i = 1}^{K} x \varphi^{(l)}_i ,\]
where  $\varphi^{(l)}_i$ is a random automorphism of  length $l$.

We use the data sets $S_P(r)$ to verify, using independent
experiments, the conclusions of Section \ref{sec:WC1} on the
average expected time $E(r)$  required for GWA to solve the length
reduction problem in the group $F_r$. If they  are true then the
expected number of iterations $Gen_r(w)$ required for GWA to
produce $w_{min}$ for a given input $w \in F_r$ satisfies the
following estimate:
 \begin{equation} \label{eq:check}
Gen_r(w) \leq E(r)CW(w) \leq E(r)|w|\end{equation}

 Let $Q_r$ be the fraction of such elements $w$ in the set $S_P(r)$ for which
  $Gen_r(w)  \leq E(r)|w|$ holds. Table \ref{tab:Epercent}  shows values of
  $Q_r$ for  $r = 2,5,10,20.$
 We can see that $Q_r$ is closed to 1 for all tested ranks, as predicted.

 In particular, we can make the following
\begin{conclusion}
\label{con:5} The genetic algorithm GWA with the termination condition
T3 gives reliable results.
\end{conclusion}

%

\begin{table}[ht]
\begin{center}
\begin{tabular}{|l|c|c|c|c|c|c|c|c|c|}
\hline
            &  $F_2$ & $F_5$ & $F_{10}$ &  $F_{15}$ & $F_{20}$\\
\hline
$E(r)$ & 1  & 3 &  7 & 12  &  18 \\
\hline
all words   & 0.93 & 0.93  & 0.99  & 0.99   & 0.99\\
\hline
$|w| > 100$ & 1.0  & 0.99 & 0.99  &  0.99 & 1.0\\
\hline

\end{tabular}
\end{center}
\caption{Fraction of elements $w \in S_P(r)$ with  $TGen_r(w) \leq
E(r)|w|$.} \label{tab:Epercent}
\end{table}
%

In constructing the set $S_P(r)$ we  select $K$  to ensure the
representativeness of characteristic $Q_r$ (see table
\ref{tab:F10test}).
%

\begin{table}[ht]
\begin{center}
\begin{tabular}{|l|c|c|c|c|c|}
\hline
 $K$    &    $Q_2$ & $Q_5$ & $Q_{10}$ & $Q_{15}$ & $Q_{20}$  \\
\hline
100 &  0.932 & 0.923 & 0.996 & 0.995 & 0.992  \\
\hline
200   & 0.93 & 0.926 & 0.996 &  0.995 & 0.993 \\
\hline
300   &  0.928 & 0.929 & 0.996 & 0.995 & 0.993\\
\hline
400 &   0.928 & 0.928 & 0.996 & 0.995 &  0.993 \\
\hline
500 &   0.93 & 0.926 & 0.996 & 0.995 & 0.993 \\
\hline
\end{tabular}
\end{center}
\caption{Values of $Q_r$ computed with different values  of $K$. }
 \label{tab:F10test}
\end{table}

The data stabilizes at $K$ = 500 and this is  the value of $K$
used in our experiments.

\section{Time complexity of GWA} \label{subsec:GWAcomp}

 It is not easy  to estimate, or even to define,
time complexity of GWA because of its stochastic nature. However,
one can estimate the time complexity of the major components  of
GWA on each given iteration.  Afterward, one may define a time
complexity function $T_{GWA}(s)$ as an average number of
iterations required by GWA to find a solution starting on a given
input $s$.

Let GWA starts to work on an input $w \in F$. Below we give some
estimates for the time required for GWA to make one iteration. It
is easy to see that the total execution time $T_{CMR}(P)$ of
Crossover, Mutation, and Replacement operators, needed to generate
the a population $P_{new}$ from a given population $P$, does not
depend on the length of the input $w$ and depends only on the
cardinality  of the population $P$ (which is fixed), and the
length $|\mu|$ of members $\mu$ of the current population $P$
(here $|\mu|$ is the length of the sequence $\mu$). Therefore, for
some constant $C_{CMR}$ the following estimate holds
$$T_{CMR}(P) \leq C_{CMR} \cdot M_P$$
where $M_P = \max \{|\mu| \mid \mu \in P\}$.

To compute $Fit(\mu)$ for  a given $\mu \in P$ it requires to run
the substitution routine SR on the input $w\mu$. Since $|wt| \leq
3|w|$ for any restricted Whitehead automorphism $t \in T$  one has
$|w\mu| \leq 3^{|\mu|}|w|$ for each $\mu \in P$. Hence the
execution time $T_{Fit}$ required to compute $Fit(\mu)$ can be
bounded from above by
\[T_{Fit} \leq C_{Fit}\cdot |w\mu| \leq C_{Fit}\cdot 3^{M_P}\cdot |w|\]

This argument shows that the time $T_{gen}(P)$ required for GWA to
generate a new population from a given one $P$ can be estimated
from above by
$$T_{gen}(P) \leq T_{CMR}(P) + T_{Fit} \leq C_{CMR} \cdot M_P +
C_{Fit} \cdot 3^{M_P} \cdot|w|.$$
In fact, the estimate $|wt| \leq 3|w|$ is very crude, as we have
seen in Section \ref{subsec:complexity} one has on  average $|wt|
\leq c_r|w|$ and the values of $c_r$ are much smaller than 3 (see
Table \ref{tab:CW-W}). So on average one can make the following
estimate:
$$T_{gen}(P)\leq C_{CMR} \cdot M_P + C_{Fit} \cdot {c_r}^{M_P} \cdot|w|.$$
Thus, the length of members  of the current population $P$ has
crucial impact on the time complexity of the procedure that
generates  the next population.

\begin{figure}[ht]
\centerline{  \includegraphics[scale=0.5]{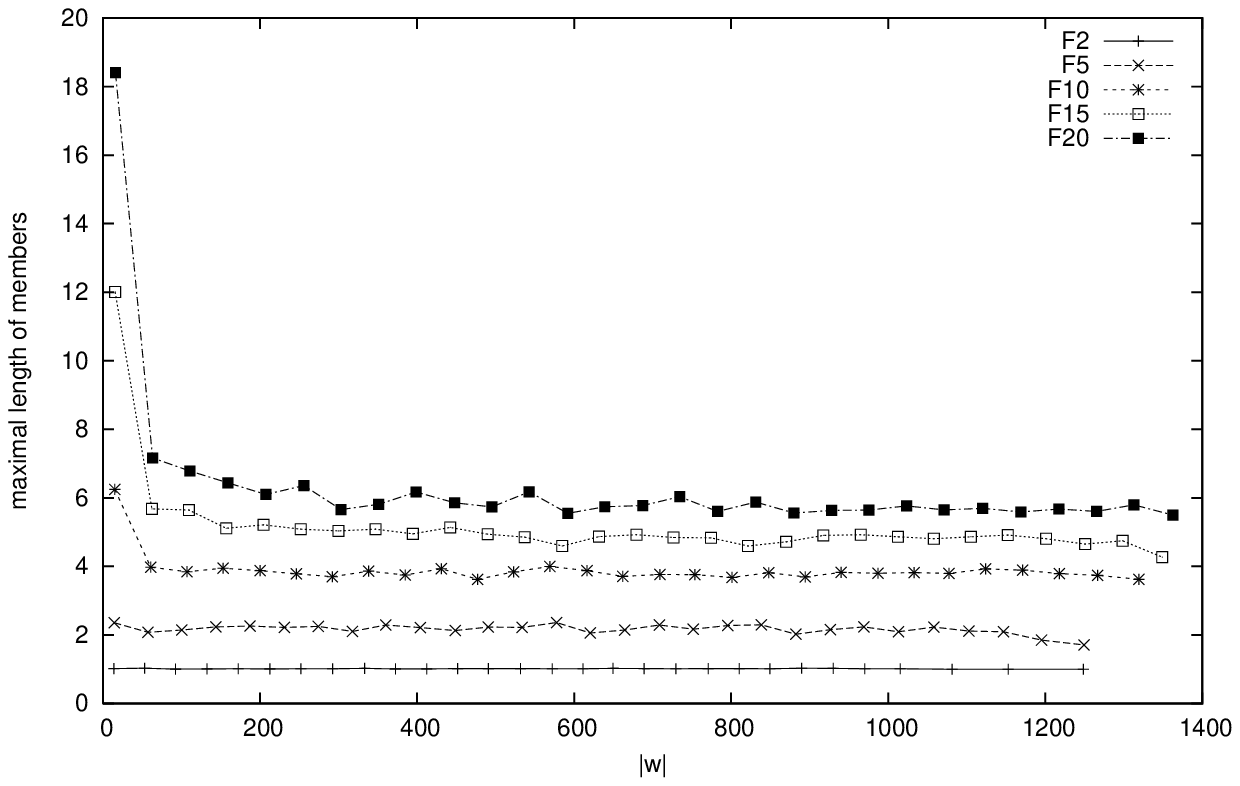}
\includegraphics[scale=0.5]{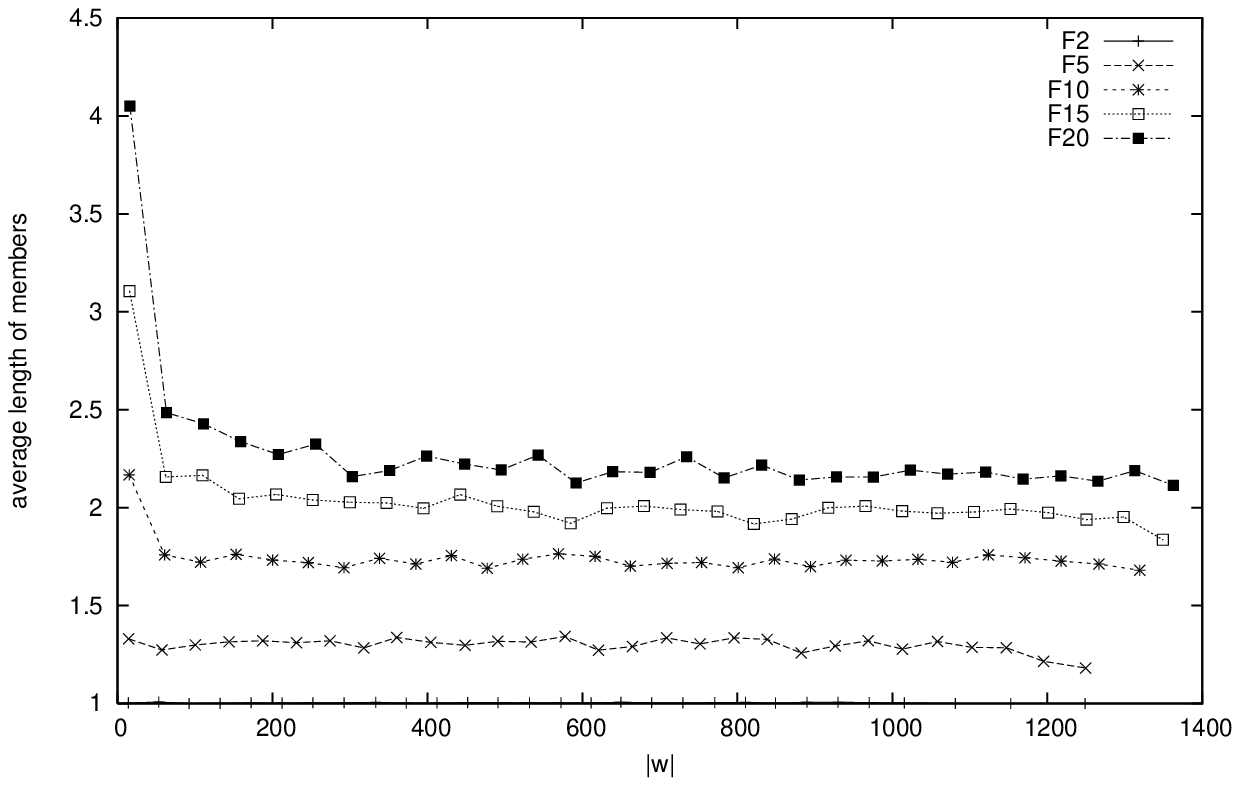}} \centerline{a)
\hspace{5.5cm}  b)}
\caption{Values of $|\mu|$ for various word lengths:
a) maximal $|\mu|$, b) average $|\mu|$.} \label{fig:mu}
\end{figure}

A priori, there are no limits on the length of the population
members $\mu \in P$. However, application of  the Substitution
Method (Section \ref{sub:substitution})   divides GWA  into  a
sequence of separate runs, each of which solves the Length
Reduction Problem for a current word $w_i = wt_1 \ldots t_i$.
Furthermore, our experiments show  that to solve this  problem
 GWA generates population members in $P$ of the  average length
 $E|\mu|$ which {\em does not depend} on the length of the input $w_i$,
 it depends only on the rank of
 $F$. In  Figure \ref{fig:mu}  we present results of our experiments
with computing $|\mu|, (\mu \in P)$ when running GWA on inputs $w$
from $S_{NMin}(r)$.

\noindent In Table \ref{tab:avg_len_m} we collect average and
maximal values of $|\mu|$ for inputs $w \in S_{NMin}(r)$ for
various ranks $r$.

%
\begin{table}[ht]
\begin{center}
\begin{tabular}{|l|c|c|c|c|c|}
\hline
                        & $F_2$   & $F_5$  &$F_{10}$ &$F_{15}$ &$F_{20}$ \\
\hline
Average $|\mu|$                   &  1.0 & 1.3 & 1.7  & 2.0  & 2.3   \\
\hline
Maximal $|\mu|$                   &  1.0 & 2.2 & 3.8  & 5.1  & 6.3   \\
\hline
\end{tabular}
\end{center}
\caption{Maximal and average lengths of the population members.}
\label{tab:avg_len_m}
\end{table}

\noindent This experimental data allows us to state   the
following observed phenomena.

\begin{conclusion}
\label{con:CompGWA} To solve the Length Reduction problem for a
given non-minimal $w \in F$ GWA  generates  new populations in
time  bounded from above by $C_r|w|$ where $C_r$ is a constant
bounded from above in the worst case by
 $$C_r \leq C_{CMR} \cdot M_P + C_{Fit} \cdot 3^{M_P},$$
 and on average by
$$C_r \leq C_{CMR} \cdot M_P + C_{Fit} \cdot c_r^{M_P},$$
\end{conclusion}

 Now we can estimate the expected  time-complexity $TGWA_r(w)$ of GWA on
 an input $w \in F_r$ as
follows:
\[TGWA_r(w) \approx Gen_r(w)\cdot average(T_{gen}(P)) \leq
E(r)\cdot WC_T(w)\cdot C_r \cdot |w|.\]

We conclude this section with a comment that average values of
$|\mu| (\mu \in P)$ shed some light on the average height of
"picks" (see Section \ref{se:problems}) for the set $T$ of
restricted Whitehead automorphisms. This topic needs a separate
research and we plan to address this issue in the future.

\subsection{Comparison of the standard Whitehead algorithm with
the genetic Whitehead algorithm} \label{se:comparison}

 In this section
we compare results of our experiments with the standard Whitehead
algorithm DWA and  the genetic algorithm GWA. We tested these
algorithms on the set $S_P$ of pseudo-random primitive elements.

As we have seen in Section \ref{subsec:GWAcomp} we may estimate
the expected  time required for GWA to find a length reducing
automorphism on a non-minimal input
 $w \in F_r$ as:
 $$C_r \cdot E(r) \cdot |w|.$$
 Recall from Section \ref{sub-sec:comp-DWA} that the expected  time
 required for  DWA to find such an automorphism
  can be estimated by
$$ \frac{A_r}{|LR_r|} \cdot |w|.$$
 In  Table \ref{tab:S_non_min_repres}
 and Figure \ref{fig:LR} we collected an
experimental data on average values of $E(r)$ and
$\frac{A_r}{|LR_r|}$ for various free groups $F_r$. It seems from
our experiments  that
 $$C_r \cdot E(r)<<  \frac{A_r}{|LR_r|}$$
 for big enough $r$. Thus, we should expect much better
performance of GWA than DWA on groups of higher ranks.

\begin{figure}[ht]
\centerline{  \includegraphics[scale=0.5]{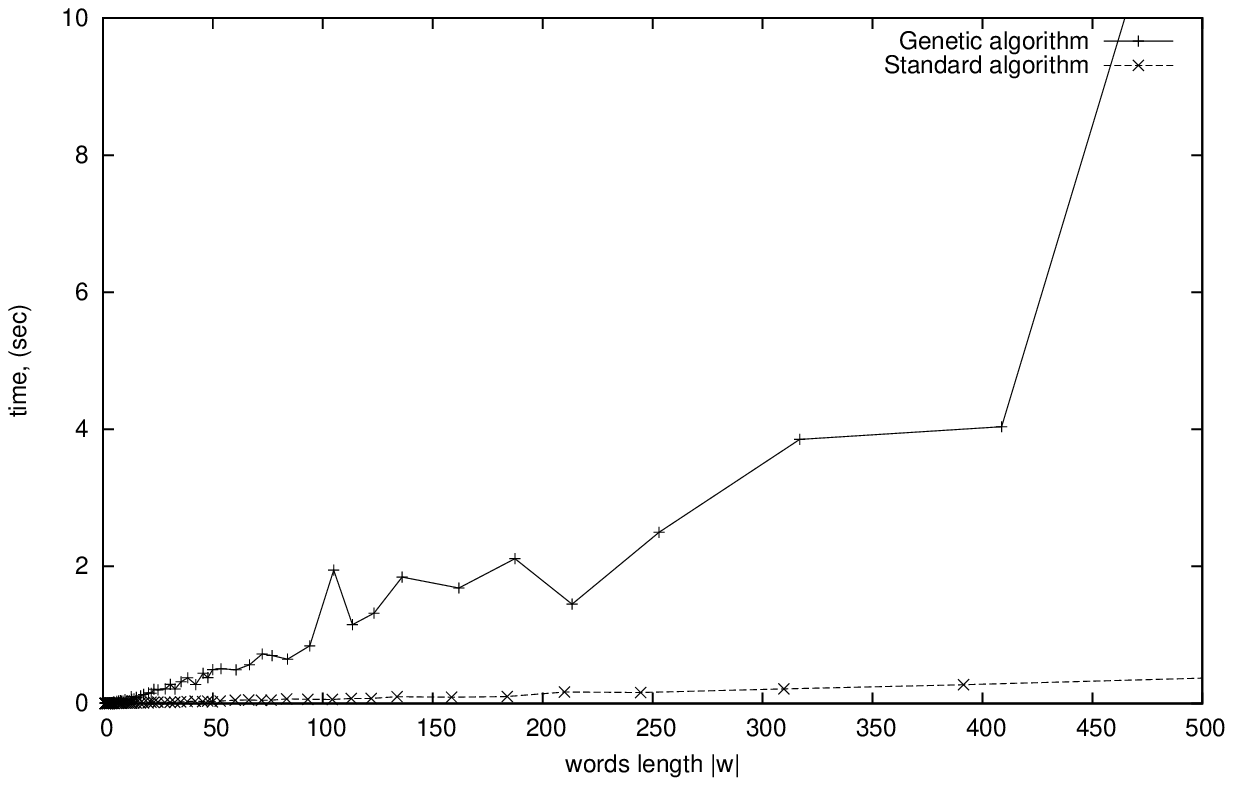}
\includegraphics[scale=0.5]{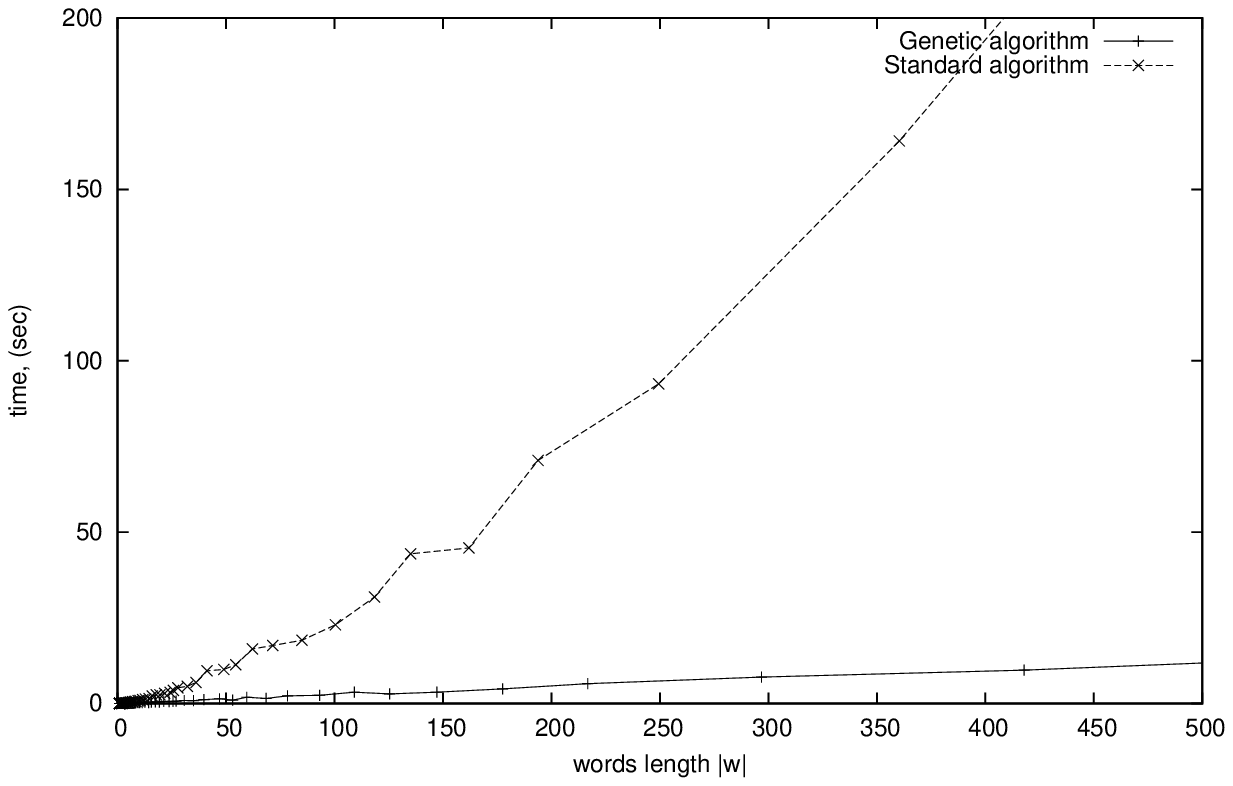}} \centerline{a)
\hspace{5.5cm}  b)}


\centerline{ \includegraphics[scale=0.5]{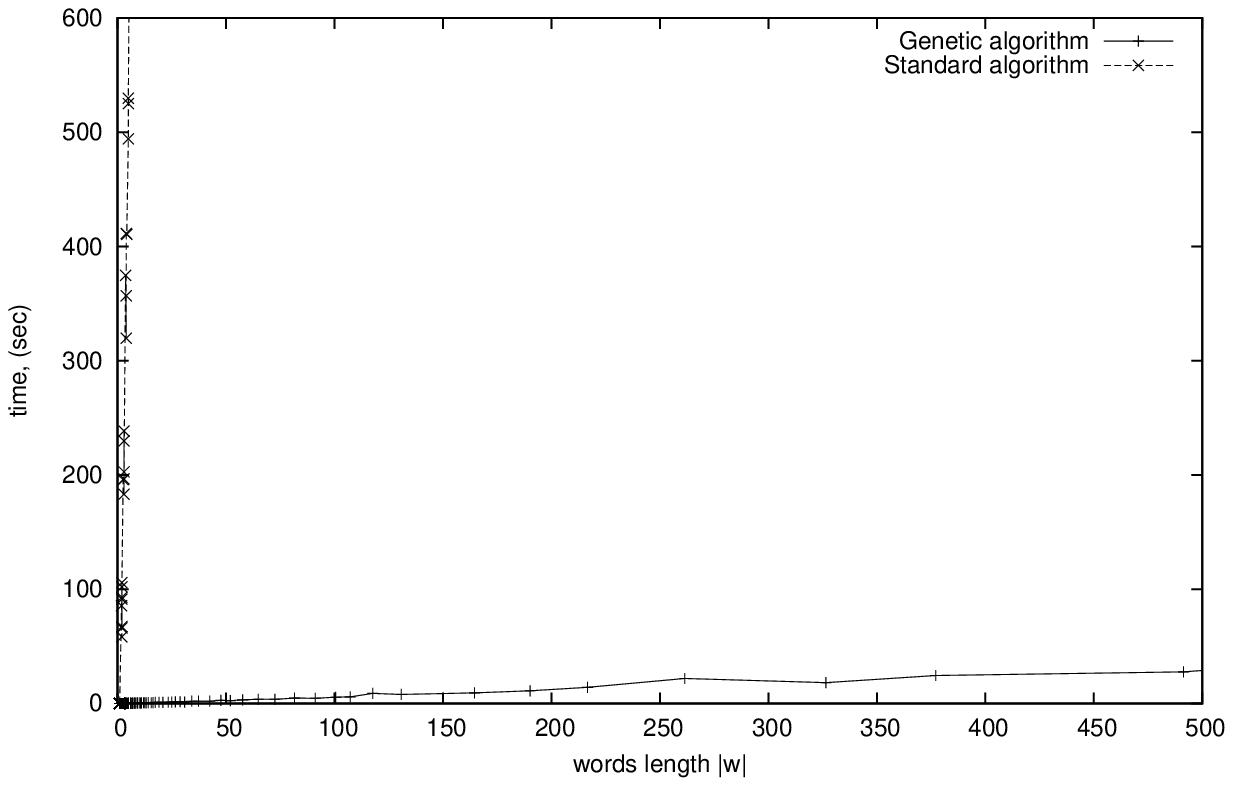}}
\centerline{ c) } \caption{Time comparison between standard and genetic
algorithms on primitive elements in a) $F_2$, b) $F_5$ and c)
$F_{10}$.} \label{fig:cmpF}
\end{figure}

In Table \ref{tab:table1} and Figures \ref{fig:cmpF} we present
results on performance comparison of  GWA with an implementation
of the standard Whitehead's algorithm DWA available in
\emph{MAGNUS} software package \cite{magnus}. We run the
algorithms on words $w \in S_P(r)$ and  measured the  execution
time.  We terminated an algorithm if it was unable to obtain the
minimal element (of length 1) on an input $w$ after being running
for more then an hour. There were very few runs of DWA for words
$w \in F_{10}$ with $|w| > 100$ that finished within an hour.
There were no such runs for $|w|
> 200$ at all, and therefore results of these experiments are marked
``na'' (not available).

%

\begin{table}[ht]
\begin{center}
\begin{tabular}{|l|c|c|c|c|c|c|c|c|c|}
\hline

& \multicolumn{3}{c|}{$F_2$} & \multicolumn{3}{c|}{$F_5$} &
\multicolumn{3}{c}{$F_{10}$} \vline \\

\hline
$|U|$            &  57 & 104 & 268&  57 & 106 & 228 & 52  & 102 & 268\\
\hline
Time spent        &     &     &    &     &     &     &     &     &\\
by the standard   & 0.03 & 0.07 & 0.18 & 13.29  & 27.4 & 85.9 &  1995& na\footnotemark &  na\\
algorithm, s      &     &     &    &     &     &     &     &     &\\
\hline
Time spent        &     &     &    &     &     &     &     &     & \\
by the genetic    & 0.52 & 1.2 & 2.7& 1.4 & 2.6& 5.6& 2.6 & 6.07  & 17.4\\
algorithm, s      &     &     &    &     &     &     &     &     &\\
\hline
\end{tabular}
\end{center}
\caption{Performance comparison of  DWA and GWA.} \label{tab:table1}
\end{table}
%


\begin{conclusion}
\label{con:6} GWA performs much better than DWA in free groups
$F_r$ for sufficiently big $r$ (in our experiments, $r \geq 5$)
and on  sufficiently long inputs (in our experiments, $|w| \geq
10$).
\end{conclusion}

\section{Mathematical problems arose from the experiments}
\label{se:problems}

 We believe that there must be some hidden
mathematical reasons for the genetic algorithm GWA to perform so
fast.  In this section we formulate several mathematical questions
which, if confirmed, would explain the robust performance of GWA,
and  lead to improved versions of the standard GWA, or to
essentially new algorithms. We focus mostly on particular choices
of the finite set of initial elementary automorphisms, and
geometry of connected components of the Whitehead graph
$\Gamma(F_r,1,\Omega_r)$.

 \begin{conjecture} \label{conj:1}
 Let $U \in F_r^k$. Then there exists a polynomial $P_{r,k}$ such that
  $$|Orb_{min}(U)| \leq P_{r,k}(|U_{min}|)$$
 \end{conjecture}

\begin{conjecture}
Almost all elements in $F_r, r \geq 2$ are Whitehead minimal.
\end{conjecture}
Of course, a rigorous formulation  of this conjecture has to
involve some   probability measure on the free group $F$. One of
the typical approaches to such problems is based on an asymptotic
density on $F$ as a measuring tool. Recently, a theoretical
justification of this conjecture, relative to the asymptotic
density, appeared in  \cite{KSS}. Below we use the asymptotic
density as our standard measuring  tool, though the measures
$\mu_s$ from \cite{BMR} would provide  more precise results.

The first conjecture deals with the average complexity of the
standard Whitehead's descent algorithm DWA.
\begin{conjecture}
\label{conj:2} Let $F = F_n$ be  a free group of rank $n$, $NMin_l
\subset F$ the set of all non-minimal elements in $F$ of length
$l$. Then there is a constant $LR_n$ such that
$$\limsup_{l \rightarrow \infty} \frac{1}{|NMin_l|}\sum_{w \in NMin_l}
|LR(w)| \ \ \ = \ LR_n.$$
\end{conjecture}

\begin{conjecture}\label{con:CW|w|}
Let
$$W_m = \{w \in F_r \mid WC(w) = m\}$$
and
$$W_{m,c_r} = \{w \in W_m \mid |w| \geq c_r^m\}$$
 There exists a constant $c_r > 1$ such that
$$\lim_{m \rightarrow \infty} \frac{|W_{m,c_r}|}{|W_m|} =
1$$
Moreover, the convergence is exponentially fast.
\end{conjecture}

Let $T = T_r$ be the  restricted set of  Whitehead automorphisms
of the group $F_r$  defined in Section \ref{sec:Sol}. Recall that
\[|T| = 5r^2 - 4r.\]
We say that $u \in Orb(w)$ is a {\em local minimum} (with respect
to the length function), if for $u \neq w_{min}$ but $|ut| \geq
|u|$ for any $t \in T$. If $u$ is a local minimum in $Orb(w)$ then
a sequence of moves $t_1, \ldots, t_k$ such that $|ut_1 \ldots
t_k| < |u|$ and $k$ is minimal with this property is called a {\em
pick} at $u$. We say that the Whitehead's descent algorithm with
respect to $T$ (see Section \ref{sub-sec:WA}) is {\em monotone} on
$w$ if it does not encounter any local minima.

\begin{conjecture}
\label{con:T-monotone}
 For "most" non-minimal elements $w \in F_r$
the Whitehead's descent algorithm with respect to $T$  is
monotone. More precisely, let $NMin_l \subset F_r$ be the set of
all non-minimal elements in $F_r$ of length $l$, and  $NMin_{l,T}$
is the subset of those for which the Whitehead's descent algorithm
with respect to $T$ is  monotone. Then
$$\lim_{m \rightarrow \infty} \frac{|NMin_{l,T}|}{|NMin_l |} = 1$$
Moreover, the convergence is exponentially fast.
\end{conjecture}

Observe, that if Conjecture \ref{con:T-monotone} holds then on
most inputs $w \in NMin \subset F_r$ the Whitehead's descent
algorithm with respect to $T$ requires at most $C\cdot r^2 \cdot
WC(w)\cdot |w|$ steps to find $w_{min}$.

 Now we are in a position to formulate  the following
conjecture
\begin{conjecture} \label{con:WCcomp-DWA}
The  time complexity  (or, at least, the average-case time
complexity) of the Problem A on inputs $w \in NMin \subset F_r$ is
bounded from above by
 $$ P(r)WC(w)|w|$$
where $P(r)$ is a fixed polynomial.
\end{conjecture}

\begin{problem}
\label{pro:1} What is geometry of the graph
$\Gamma(F_r,1,\Omega_r)$? In particular, are connected components
of  $\Gamma(F_r,1,\Omega_r)$ hyperbolic?
\end{problem}

If uncovered, the geometric properties of the graphs
$\Gamma(F_r,1,\Omega_r)$ should provide fast deterministic
algorithms for Problems A and B.

\end{document}